\documentclass[12pt]{amsart}
\usepackage{amsmath,amssymb,amsbsy,amsfonts,latexsym,amsopn,amstext,
                                               amsxtra,euscript,amscd,bm}
\usepackage{url}
\usepackage[colorlinks,linkcolor=blue,anchorcolor=blue,citecolor=blue]{hyperref}
\usepackage{color}
\usepackage{float}

\makeatletter

\newcommand{\Rmnum}[1]{\expandafter\@slowromancap\romannumeral #1@}
\makeatother

\hypersetup{breaklinks=true}

\begin{document}

\newtheorem{theorem}{Theorem}
\newtheorem{lemma}[theorem]{Lemma}
\newtheorem{claim}[theorem]{Claim}
\newtheorem{cor}[theorem]{Corollary}
\newtheorem{prop}[theorem]{Proposition}
\newtheorem{definition}{Definition}
\newtheorem{question}[theorem]{Question}
\newtheorem{remark}{Remark}
\newcommand{\hh}{{{\mathrm h}}}

\numberwithin{equation}{section}
\numberwithin{theorem}{section}
\numberwithin{remark}{section}
\numberwithin{table}{section}

\def\sssum{\mathop{\sum\!\sum\!\sum}}
\def\ssum{\mathop{\sum\ldots \sum}}
\def\iint{\mathop{\int\ldots \int}}

\def\squareforqed{\hbox{\rlap{$\sqcap$}$\sqcup$}}
\def\qed{\ifmmode\squareforqed\else{\unskip\nobreak\hfil
\penalty50\hskip1em\null\nobreak\hfil\squareforqed
\parfillskip=0pt\finalhyphendemerits=0\endgraf}\fi}

\newfont{\teneufm}{eufm10}
\newfont{\seveneufm}{eufm7}
\newfont{\fiveeufm}{eufm5}
%
%
\newfam\eufmfam
     \textfont\eufmfam=\teneufm
\scriptfont\eufmfam=\seveneufm
     \scriptscriptfont\eufmfam=\fiveeufm
%
%
\def\frak#1{{\fam\eufmfam\relax#1}}

\newcommand{\bflambda}{{\boldsymbol{\lambda}}}
\newcommand{\bfmu}{{\boldsymbol{\mu}}}
\newcommand{\bfxi}{{\boldsymbol{\xi}}}
\newcommand{\bfrho}{{\boldsymbol{\rho}}}

\def\fK{Frak K}
\def\fT{Frak{T}}

\def\fA{{Frak A}}
\def\fB{{Frak B}}
\def\fC{{Frak C}}

\def \balpha{\bm{\alpha}}
\def \bbeta{\bm{\beta}}
\def \bgamma{\bm{\gamma}}
\def \blambda{\bm{\lambda}}
\def \bchi{\bm{\chi}}
\def \bphi{\bm{\varphi}}
\def \bpsi{\bm{\psi}}

\def\eqref#1{(\ref{#1})}

\def\vec#1{\mathbf{#1}}



\def\sA{\mathscr A}
\def\sB{\mathscr B}
\def\sC{\mathscr C}
\def\sD{\mathscr D}
\def\sE{\mathscr E}
\def\sF{\mathscr F}
\def\sG{\mathscr G}
\def\sH{\mathscr H}
\def\sI{\mathscr I}
\def\sJ{\mathscr J}
\def\sK{\mathscr K}
\def\sL{\mathscr L}
\def\sM{\mathscr M}
\def\sN{\mathscr N}
\def\sO{\mathscr O}
\def\sP{\mathscr P}
\def\sQ{\mathscr Q}
\def\sR{\mathscr R}
\def\sS{\mathscr S}
\def\sU{\mathscr U}
\def\sT{\mathscr T}
\def\sV{\mathscr V}
\def\sW{\mathscr W}
\def\sX{\mathscr X}
\def\sY{\mathscr Y}
\def\sZ{\mathscr Z}


\def\fA{\mathfrak A}
\def\fB{\mathfrak B}
\def\fC{\mathfrak C}
\def\fD{\mathfrak D}
\def\fE{\mathfrak E}
\def\fF{\mathfrak F}
\def\fG{\mathfrak G}
\def\fH{\mathfrak H}
\def\fI{\mathfrak I}
\def\fJ{\mathfrak J}
\def\fK{\mathfrak K}
\def\fL{\mathfrak L}
\def\fM{\mathfrak M}
\def\fN{\mathfrak N}
\def\fO{\mathfrak O}
\def\fP{\mathfrak P}
\def\fQ{\mathfrak Q}
\def\fR{\mathfrak R}
\def\fS{\mathfrak S}
\def\fU{\mathfrak U}
\def\fT{\mathfrak T}
\def\fV{\mathfrak V}
\def\fW{\mathfrak W}
\def\fX{\mathfrak X}
\def\fY{\mathfrak Y}
\def\fZ{\mathfrak Z}


\def\C{{\mathbb C}}
\def\F{{\mathbb F}}
\def\K{{\mathbb K}}
\def\L{{\mathbb L}}
\def\N{{\mathbb N}}
\def\P{{\mathbb P}}
\def\Nc{\N^{\hskip1pt c}}
\def\Q{{\mathbb Q}}
\def\R{{\mathbb R}}
\def\Z{{\mathbb Z}}

\def\eps{\varepsilon}
\def\mand{\qquad\mbox{and}\qquad}
\def\mor{\qquad\mbox{or}\qquad}
\def\\{\cr}
\def\({\left(}
\def\){\right)}
\def\[{\left[}
\def\]{\right]}
\def\<{\langle}
\def\>{\rangle}
\def\fl#1{\left\lfloor#1\right\rfloor}
\def\rf#1{\left\lceil#1\right\rceil}
\def\le{\leqslant}
\def\ge{\geqslant}
\def\ds{\displaystyle}
\def\lealmost{\preccurlyeq}
\def\gealmost{\succcurlyeq}

\def \LLN {\prec \hskip-6pt \prec}

\def\cg{{\mathcal g}}

\def\vr{\mathbf r}

\def\e{{\mathbf{\,e}}}
\def\ep{{\mathbf{\,e}}_p}
\def\em{{\mathbf{\,e}}_m}

\def\Tr{{\mathrm{Tr}}}
\def\Nm{{\mathrm{Nm}}}

 \def\SS{{\mathbf{S}}}

\def\lcm{{\mathrm{lcm}}}

\def\fl#1{\left\lfloor#1\right\rfloor}
\def\rf#1{\left\lceil#1\right\rceil}

\def\mand{\qquad \mbox{and} \qquad}

\newcommand{\commK}[1]{\marginpar{%
\begin{color}{red}
\vskip-\baselineskip 
\raggedright\footnotesize
\itshape\hrule \smallskip K: #1\par\smallskip\hrule\end{color}}}

\newcommand{\commI}[1]{\marginpar{%
\begin{color}{magenta}
\vskip-\baselineskip 
\raggedright\footnotesize
\itshape\hrule \smallskip I: #1\par\smallskip\hrule\end{color}}}

\newcommand{\commT}[1]{\marginpar{%
\begin{color}{blue}
\vskip-\baselineskip 
\raggedright\footnotesize
\itshape\hrule \smallskip T: #1\par\smallskip\hrule\end{color}}}




\hyphenation{re-pub-lished}

\mathsurround=1pt

\def\bfdefault{b}
\overfullrule=5pt

\def \F{{\mathbb F}}
\def \K{{\mathbb K}}
\def \Z{{\mathbb Z}}
\def \Q{{\mathbb Q}}
\def \R{{\mathbb R}}
\def \C{{\\mathbb C}}
\def\Fp{\F_p}
\def \fp{\Fp^*}

\def\Kmn{\cK_p(m,n)}
\def\psmn{\psi_p(m,n)}
\def\SI{\cS_p(\cI)}
\def\SIJ{\cS_p(\cI,\cJ)}
\def\SAIJ{\cS_p(\cA;\cI,\cJ)}
\def\SABIJ{\cS_p(\cA,\cB;\cI,\cJ)}
\def \xbar{\overline x_p}

\title[Squares in Piatetski--Shapiro sequences]{Squares in Piatetski--Shapiro Sequences}

 \author[K. Liu] {Kui Liu}
\address{School of Mathematics and  Statistics, Qingdao University, No.308, Ningxia Road, Shinan, Qingdao, Shandong, 266071, P. R. China}
\email{liukui@qdu.edu.cn}

 \author[I. E. Shparlinski] {Igor E. Shparlinski}
\thanks{T. P. Zhang is the corresponding author (tpzhang@snnu.edu.cn).}

\address{Department of Pure Mathematics, University of New South Wales,
Sydney, NSW 2052, Australia}
\email{igor.shparlinski@unsw.edu.au}

 \author[T. P. Zhang] {Tianping Zhang}
\address{School of Mathematics and Information Science, Shaanxi Normal University, Xi'an 710119 Shaanxi, P. R. China}
\email{tpzhang@snnu.edu.cn}

\begin{abstract} We study
the distribution of squares in a Piatetski-Shapiro sequence
$\(\fl{n^c}\)_{n\in\N}$  with $c>1$ and $c\not\in\N$. We also study more general
equations $\fl{n^c} = sm^2$, $n,m\in \N$, $1\le n \le N$ for an integer $s$ and
obtain several bounds on the number of solutions for a fixed $s$ and on average
over $s$ in an interval.  These results are based on various techniques chosen depending
on the range of the parameters.
\end{abstract}


\maketitle

\section{Introduction}

\subsection{Motivation and formulation of the problem}

{\it Piatetski-Sh\-apiro sequ\-en\-ces\/} (PS-sequ\-en\-ces), that is,  sequences of the form
$$
\N^c=(\fl{n^c})_{n\in\N}\qquad (c>1,\ c\not\in\N),
$$
where $\fl{z}$ is the integer part of a real $z$,
have been extensively studied by many authors since their introduction by Piatetski-Shapiro~\cite{PS},
see~\cite{AkGu,Bak,BaBa,BBBSW,BBGY,BGS,Guo,Spi} and the references therein.

Here we consider the distribution of perfect squares in PS-sequences,
which seems to be a new, yet natural question to study.  More precisely,
for a real $c>1$ and positive integers $N$ and $s$,  we denote by
$Q_{c}(s;N)$ the number of solutions to the equation
$$
\fl{n^c} = sm^2, \qquad  1 \le n \le N, \ m,n \in \Z.
$$
Clearly, we have the following trivial bound
\begin{equation}
\label{eq:Triv}
Q_{c}(s;N)  \le \min\left\{N,\ s^{-1/2}N^{c/2}\right\}.
\end{equation}
Here we use a variety of different techniques to obtain asymptotic formulas,
or upper bounds improving~\eqref{eq:Triv}. We also  study $Q_{c}(s;N) $
on average over positive  square-free integers $s \le S$, that is, the quantity
$$
\fQ_c(S,N) = \sum_{\substack{s\leq S\\s~\text{square-free}}} Q_{c}(s;N).
$$
We remark that only the case $S\leq N^c$ is meaningful, hence we always assume this.
Having nontrivial upper bounds on $\fQ_c(S,N)$ immediately implies a lower
bound on the number of distinct square-free parts of the integers $\fl{n^c}$,
$1 \le n \le N$. In turn, this can be reformulated as a lower bound on the
number of distinct quadratic fields in the sequence of fields $\Q\(\sqrt{\fl{n^c}}\)$,
$1 \le n \le N$.

\subsection{Notation}
Among other methods, our results are also based on the square sieve of
Heath-Brown~\cite{HB2}
coupled with a bound of character sums with PS-sequences due to Baker and  Banks~\cite{BaBa}.
We also employ the method of exponent pairs, we refer
to~\cite[Chapter~3]{GrKol}, \cite[Sections~7.3 and~17.4]{Huxley},
\cite[Chapter~8]{IwKow} and~\cite[Chapter~3]{Mont} an exact definition, properties  and examples of exponent pairs.

Throughout the paper,  as usual $U\ll V$ and $U = O(V)$ are  both equivalent to the inequality $|U|\le B V$
with some  constant $B>0$, which maybe depend on the parameter $c$ (and sometimes, where
obvious, on the some other auxiliary parameters), however it is always uniform with respect to
our main parameters $N$, $s$ and $S$.

For two quantities $U$ and $V$, which among other parameters also depend on $N$, we use $U \LLN V$
to denote that $U \le V N^{o(1)}$  as $N\to \infty$.

We also write $u \sim U$ to denote that $U < u \le 2U$.

The letters $\ell$ and $p$, with or without subscripts, always denote prime numbers.

As usual $(k/q)$ denotes the {\it Jacobi symbol\/} modulo $q$, which we use
only for prime $q$, when  it is called the {\it Legendre symbol\/}, or for products of two primes).

We also use $\Box$ to denote a nonspecified integer square, that is,   $n=\Box$ is
equivalent to the statement that $n$ is a perfect square and thus we can write
$$ Q_c(s;N)=\sum_{\substack{n\leq N \\ \lfloor n^c\rfloor=s \Box}}1.$$

\section{Main Results}

\subsection{Results for a fixed $s$}

We start with an asymptotic formula for $Q_{c}(s;N)$ for the values of $c$ close to $1$.
We refer to~\cite[Chapter~3]{GrKol}, \cite[Sections~7.3 and~17.4]{Huxley},
\cite[Chapter~8]{IwKow} and~\cite[Chapter~3]{Mont}  for a background on
exponent pairs.

\begin{theorem}
\label{thm: Asymp s}
For any $c > 1$, $c\not\in\N$ and any exponent pair  $(\kappa,\lambda)$
we have an asymptotic formula
\begin{equation*}
\begin{split}
Q_c(s;N)=\ &\gamma(2\gamma-1)^{-1}s^{-1/2}N^{1-c/2}\\
& +O\(s^{-\rho_1(\kappa, \lambda)}N^{\vartheta_1 (c,\kappa,\lambda)+o(1)}+s^{-\rho_2(\kappa, \lambda)}N^{\vartheta_2 (c,\kappa,\lambda)+o(1)}\)
\end{split}
\end{equation*}
as $N \to \infty$, where $\gamma=1/c$,
$$
\rho_1(\kappa, \lambda) = \frac{\lambda}{2(\kappa+1)},\indent
\vartheta_1 (c,\kappa,\lambda) = \frac{2\kappa + c \lambda}{2(1+\kappa)}
$$
and
$$
\rho_2(\kappa, \lambda) = \frac{(\lambda-\kappa)}{2},\indent
\vartheta_2 (c,\kappa,\lambda) = \frac{2\kappa+c(\lambda-\kappa)}{2}.
$$
\end{theorem}

For example, taking $(\kappa,\lambda)=(9/56, 37/56)$ (see~\cite[Chapter~7]{GrKol}), we have
\begin{equation*}
\begin{split}
Q_c(s;N)=\ &\gamma(2\gamma-1)^{-1}s^{-1/2}N^{1-c/2}\\
& + O\(s^{-37/130}N^{(18 + 37 c)/130+o(1)}+s^{-1/4}N^{(9+14c)/56+o(1)}\),
\end{split}
\end{equation*}
which gives an asymptotic formula for $s=1$ and $1<c<56/51\approx 1.09804$.

Furthermore, taking $(\kappa,\lambda)=(1/2,1/2)$  (see~\cite[Chapter~3]{GrKol}),
we obtain
\begin{equation*}
\begin{split}
Q_c(s;N)&=\gamma(2\gamma-1)^{-1} s^{-1/2}N^{1-c/2} \\
&\qquad \qquad \qquad  +O\(s^{-1/6}N^{(2+c)/6+o(1)}+N^{1/2+o(1)}\)\\
&\LLN s^{-1/6}N^{(2+c)/6}+N^{1/2}
\end{split}
\end{equation*}
for $c>1$,  $c\notin\N$ (which is nontrivial for $s=1$ and $1< c < 4$, $c\notin\N$).

%
%

For larger values of $c$ we have a less explicit bound, which is
nontrivial for any $c>2$.
This bound depends on an absolute constant $\beta(c)>0$, depending only on $c$ such that
for any positive  integers $N$ and $q$, for characters sums
\begin{equation}
\label{eq:sum T}
T_{c,\chi}(q;N)=\sum_{N/2 < n\le N}\chi(\lfloor n^c\rfloor),
\end{equation}
with  a primitive Dirichlet character  $\chi$ modulo $q$
(see~\cite[Chapter~3]{IwKow} for a background on characters),
we have
\begin{equation}
\label{eq:Betac}
T_{c,\chi}(q;N)\ll q^{1/2}N^{1-\beta(c)}.
\end{equation}
The existence of such $\beta(c)$  for any $c>2$ of the form
\begin{equation}
\label{eq:beta}
\beta(c) = \beta/c^2
\end{equation}
with an absolute constant $\beta > 0$ is essentially
 a result of Baker and Banks~\cite[Theorem~1.6]{BaBa},
which we also present as Lemma~\ref{lem:q N large c} below.

\begin{theorem}
\label{thm: U-Bound}
For any   $c > 2$, $c\not\in\N$ and $\beta(c)$ satisfying~\eqref{eq:Betac}, for $N \to \infty$ we have
$$
Q_c(s; N) \ll N^{1-\beta(c)/2+o(1)}.
$$
\end{theorem}

We note that the proof of Theorem~\ref{thm: U-Bound} is based on the square sieve method of Heath-Brown~\cite{HB2} which seems to be the first application of this method
in the context of PS-sequences for large $c$, where usually the method of exponential sum
is used for small $c$. This has become possible because of the recent results of
Baker and Banks~\cite{BaBa}.

\subsection{Results on average over $s$}

Here we show that using a result of Fouvry and Iwaniec~\cite[Theorem~3]{FoIw},
when $c$ is near to $1$, we can take advantage of averaging over $s$ and estimate the sum $\fQ_c(S,N)$
better than via a direct applications of Theorem~\ref{thm: Asymp s}.

Our result, as it is natural to expect, depends on the following  function
$$
\Phi(S) = \sum_{\substack{s\leq S\\ s~\text{square-free}}}s^{-1/2}.
$$


Using the well known result, see~\cite[Theorem~333]{HaWr},
$$
\sum_{\substack{s\leq t\\ s~{\text{square-free}}}}1=\frac{6}{\pi^2} t +O(\sqrt t)
$$
and partial summation, we easily derive
\begin{equation}
\label{eq:Phi(S)}
\Phi(S)  = \frac{12}{\pi^2} S^{1/2} +O(\log S),
\end{equation}
which we can use together with the bound of Theorem~\ref{thm: Asymp s}.

\begin{theorem}
\label{thm: U-Asymp Aver}
For any $c > 1$, $c\not\in\N$, for $N \to \infty$ we have
\begin{align*}
&\fQ_c(S,N)- \frac{ 12 \gamma}{\pi^2(2\gamma-1)} S^{1/2}N^{1-c/2}\\
&\qquad\LLN S^{1/5}N^{(1+2c)/5}+S^{5/8}N^{3c/8}+S^{1/8}N^{(2+3c)/8}+SN^{1-c}
\end{align*}
with $\gamma=1/c$.
\end{theorem}

In particular, we have:
\begin{cor}
\label{cor: U-Bound Aver}
For any $c > 1$, $c\not\in\N$, any $\varepsilon>0$ and $S \leq N^{\tau(c)-\varepsilon},$
where
$$
\tau(c)= \left\{
\begin{array}{ll}
(8-3c)/5,  &{\text{for\ }} 1<c\leq 12/7,\\
2(2-c),  &{\text{for\ }} c>12/7,\\
\end{array}
\right.
$$
we have
$$
\fQ_c(S,N) = o(N)
$$
as $N \to \infty$.
\end{cor}
Clearly Corollary~\ref{cor: U-Bound Aver} is nontrivial only for $c<2$, $c\notin \N$.

\begin{remark}
 Applying Theorem~\ref{thm: Asymp s} with $(\kappa, \lambda)=(1/2, 1/2)$ and trivial estimate, we may take $\tau(c) = (4-c)/5$, which is nontrivial for $c<4$, $c\notin \N$. But for large $c$, we need refer to the square sieve again to get  a positive value for $\tau(c)$.
\end{remark}

\begin{theorem}
\label{thm: U-Bound Aver}
For any $c > 2$, $c\not\in\N$  and $\beta(c)$ satisfying~\eqref{eq:Betac},  for $N \to \infty$ we have
$$
\fQ_c(S,N) \LLN SN^{1-\beta(c)}+S^{3/4}N^{1-\beta(c)/2}.
$$
\end{theorem}

In particular, we have:

\begin{cor}
\label{cor: U-Bound Aver c>2}
For any $c > 2$, $c\not\in\N$,  $\beta(c)$ satisfying~\eqref{eq:Betac},  and any $\varepsilon>0$,
for $S \leq N^{2\beta(c)/3-\varepsilon}$,
we have
$$
\fQ_c(S,N) = o(N)
$$
as $N \to \infty$,
\end{cor}

Corollary~\ref{cor: U-Bound Aver} and Corollary~\ref{cor: U-Bound Aver c>2}
cover the full range $c>1$, $c\not\in\N$,  provided~\eqref{eq:Betac}  holds.
Hence, combining this with~\eqref{eq:beta}, which we have by Lemma~\ref{lem:q N large c} below, we obtain:

\begin{cor}
\label{cor: large square-free part}
For any $c>1$, $c\not\in\N$, there exists a constant $\vartheta(c)>0$ such that the square-free parts of almost all integers of the type $\fl{n^c},n\leq N$ are larger than $N^{\vartheta(c)}$.
\end{cor}

\section{Preparations}

\subsection{Some general statements}

As usual, we define the function $\psi(u) = u -\fl{u}-1/2$.
We use the following  result of Vaaler~\cite{Vaal},
see also~\cite[Theorem~A.6]{GrKol}.

\begin{lemma}
\label{lem:Vaal}
Let $H\geq 1$. There are functions $a(h)$ and $b(h)$,   such that for $1\leq |h|\leq H$
we have
$$a(h)\ll \frac{1}{|h|} ,\qquad  b(h)\ll \frac{1}{H},
$$
and
$$
\left|\psi(t)-\sum_{1\le |h|\leq H}a(h)e(ht)\right|\leq\sum_{|h|\leq H}b(h)e(ht).
$$
\end{lemma}

Note that we can take explicitly,
$$
a(h)=(2\pi i h)^{-1}F\(\frac{h}{H+1}\)\quad \text{and} \quad  b(h)=\frac{1}{2H+2}\(1-\frac{|h|}{H+1}\),
$$
with $F(u)=\pi u(1-|u|)\cot(\pi u)$ in  Lemma~\ref{lem:Vaal}. We also remark the right hand side
of  Lemma~\ref{lem:Vaal} is a real nonnegative number, so now absolute value symbol is necessary.
It is also important to notice that the summation on  the right hand side also includes $h=0$.

We also need the following technical result, see~\cite[Lemma~2.4]{GrKol}.

\begin{lemma}
\label{lem: Optim}
Let
$$
L(Z)=\sum_{i=1}^uA_i Z^{a_i}+\sum_{j=1}^{v}B_j Z^{-b_j},
$$
where $A_i,B_j,a_i$ and $b_j$ are positive. Let $0\leq Z_1\leq Z_2$. Then there is some $Z\in(Z_1,Z_2]$ with
$$
L(Z)\ll \sum_{i=1}^u\sum_{j=1}^v
\(A_i^{b_j}B_j^{a_i}\)^{1/(a_i+b_j)}+\sum_{i=1}^uA_iZ_1^{a_i}+\sum_{j=1}^vB_jZ_2^{-b_j},
$$
where the implied constant depends only on $u$ and $v$.
\end{lemma}

The following is a form of the square sieve  of
Heath-Brown~\cite{HB2} which is given by Friedlander
and Iwaniec~\cite[Proposition~3.1]{FrIw}
(combined with the trivial observation that if for some
integers $r$ and $s$ we have $r = s\Box$ then $rs = \Box$).
\begin{lemma}
\label{lem:Sq Sieve}
Let $a_r$,  $r=1,\ldots, R$, be an arbitrary finite sequence non-negative real numbers and let $P\geq 2$. Then we have
$$
\sum_{\substack{r=1\\ r=s\Box}}^R
a_r\leq 10P^{-2}\sum_{r=1}^R a_r\(\(\sum_{P<p\leq 2P}\(\frac{sr}{p}\)\log p\)^2+(\log sr)^2\).
$$
\end{lemma}

One can apply Lemma~\ref{lem:Sq Sieve} directly to $Q_c(s; N)$ but it is
technically easier to work with dyadic intervals, so we define
\begin{equation}
\label{eq:Q*}
 Q_c^*(s;N)=\sum_{\substack{N/2 < n\leq N \\ \lfloor n^c\rfloor=s \Box}}1.
\end{equation}

Taking $a_r$ to be the characteristic function of the event $r = \fl{n^c}$ for some positive
integer $N/2 < n \le N$ we obtain:

\begin{cor}
\label{cor:Sq Sieve PS}
For any positive integers $N$, $P\ge 2$ and  $s$,  we have
$$
Q_c^*(s; N) \ll P^{-2}\sum_{N/2 < n\leq N}\(\sum_{P<p\leq 2P}\(\frac{s \lfloor n^c\rfloor}{p}\)\log p\)^2+ P^{-2}N\log^2 N.
$$
\end{cor}

We need the following mean value estimate for real character sums,
which is Theorem 1 of ~\cite{HB3} (see also~\cite[Theorem~7.20]{IwKow}).

\begin{lemma}
\label{Heath-Brown char sum}
For any integers $M,N\ge 1$ and  complex numbers $a_n$, $n=1, \ldots, N$ we have
$$
\sum_{\substack{m\leq M\\m~\text{square-free}}}
\left|\sum_{\substack{n\leq N\\n~\text{square-free}}}a_n\left(\frac{n}{m}\right)\right|^2\le
(MN)^{o(1)}(M+N)\sum_{n\leq N}|a_n|^2,
$$
as $MN\to \infty$.
\end{lemma}

\subsection{Character sums with PS-sequences}
\label{sec:PS char sum}

We now recall the following  bound  on the sums  $T_{c,\chi}(q;x)$
defined by~\eqref{eq:sum T}, given by
Baker and Banks~\cite[Theorem~1.6]{BaBa}
(used with $y=x=N/2$), which is nontrivial for any $c>2$
(provided that $N$ is sufficiently large compared to $q$).

\begin{lemma}
\label{lem:q N large c}
Let  $N \ge 2$ and $q\geq 3$. Then for $c>2$, $c\not\in\N$, there exists an absolute
constant $\beta>0$ such that
$$
T_{c,\chi}(q;N)\ll q^{1/2}N^{1-\beta/c^2}.
$$
\end{lemma}

In particular,   Lemma~\ref{lem:q N large c}  shows
that~\eqref{eq:beta} is satisfied for some $\beta>0$ and thus the
assumption~\eqref{eq:Betac} is not void.

There is no doubt that the value of $\beta$ in Lemma~\ref{lem:q N large c}  can be explicitly evaluated.

\subsection{Exponential sums with monomials}
\label{sec:Exp sum monom}

We need the following bound due to Fouvry and Iwaniec~\cite[Theorem~3]{FoIw}.
We remark that the more recent bound of Robert and Sargos~\cite[Theorem 1]{RobSar}
does not bring any improvement to our results (as the bounds of~\cite[Theorem~3]{FoIw}
and~\cite[Theorem 1]{RobSar} have some common terms and these are exactly the
terms that dominate in our applications).

\begin{lemma}\label{Iwaniec's lemma}
Let $\alpha,\alpha_1,\alpha_2$ be real constants such that
$$\alpha\neq 1 \mand \alpha\alpha_1\alpha_2\neq 0.
$$
Let $M,M_1,M_2,x\geq 1$ and let
$$\Phi=\(\varphi_m\)_{m \sim M} \mand  \varPsi=\(\psi_{m_1,m_2}\)_{m_1\sim M_1, m_2\sim M_2}
$$
be
 two sequences of complex numbers supported on $m \sim M$, $m_1\sim M_1$ and $m_2\sim M_2$
 with $|\varphi_{m}|\leq 1$ and $|\psi_{m_1,m_2}|\leq 1$. Then for the sum
 $$
S_{\Phi, \varPsi}(x; M,M_1,M_2)=\sum\limits_{m\sim M}\sum\limits_{m_1\sim M_1}\sum\limits_{m_2\sim M_2}\varphi_m\psi_{m_1,m_2}e\left(x\frac{m^{\alpha}m_1^{\alpha_1}m_2^{\alpha_2}}{M^{\alpha}M_1^{\alpha_1}M_2^{\alpha_2}}\right),
$$
we have
 \begin{equation*}
\begin{split}
S_{\Phi, \varPsi}(x&; M,M_1,M_2)\\
& \ll\bigl(x^{1/4}M^{1/2}(M_1M_2)^{3/4}+M^{7/10}M_1M_2+M(M_1M_2)^{3/4}\\
&\qquad \qquad \qquad \qquad \quad  +x^{-1/4}M^{11/10}M_1M_2\big)\log^2(2MM_1M_2).
\end{split}
\end{equation*}
\end{lemma}

\section{Proofs of  main results}

\subsection{Proof of Theorem~\ref{thm: Asymp s}}

For $c>1$,  $c\notin\N$, let $\gamma=1/c$.
It is easy to see that $\lfloor n^c\rfloor=s\Box$ if and only if
$$
(s\Box)^{\gamma}\leq n<(s\Box+1)^{\gamma}.
$$
We also set $M = s^{-1/2}N^{c/2}$.
Using a similar argument as that in Heath-Brown~\cite{HB1}
(and in many other works on PS-sequences),
we have
$$
Q_{c}(s;N)=\sum_{m\leq M}\(\lfloor -s^\gamma m^{2\gamma}\rfloor-\lfloor -(sm^2+1)^{\gamma}\rfloor\)+O(1).
$$
Let
$$
\psi(x)=x- \fl{x}-1/2.
$$
Then we obtain
\begin{equation}
\begin{split}
\label{eq:Expand}
Q_{c}(s;N)=\sum_{m\leq M}&\((sm^2+1)^{\gamma}-s^\gamma m^{2\gamma}\)\\
&-\sum_{m\leq M}\psi(-s^\gamma m^{2\gamma})\\
& \qquad +\sum_{m\leq M}\psi\(-\(sm^{2}+1\)^\gamma\)+O(1).
\end{split}
\end{equation}
The first term in the right side is
\begin{equation}
\begin{split}
\label{eq:MainTerm}
\sum_{m\leq M}&\((sm^2+1)^{\gamma}-s^\gamma m^{2\gamma}\) \\
& =\sum_{m\leq M} s^\gamma m^{2\gamma}\(\gamma m^{-2}s^{-1} +O\(m^{-4}s^{-2}\)\)\\
&=\gamma(2\gamma-1)^{-1} s^{\gamma  -1}(N^{c/2}s^{-1/2})^{2\gamma-1}+O\left(1\right)\\
&= \gamma(2\gamma-1)^{-1} N^{1-c/2} s^{-1/2} + O\left(1\right),
\end{split}
\end{equation}
which gives the desired main term.

Now we need to estimate the other two sums with the  $\psi$-functions. We only estimate the  first sum with
$\psi\(-s^\gamma m^{2\gamma}\)$ and the other sum  with $\psi\(-\(sm^{2}+1\)^\gamma\)$ can be treated similarly and
admits the  same upper bound.

We now fix some parameter $H \ge 1$ and using Lemma~\ref{lem:Vaal}, we obtain
\begin{equation}
\label{eq:M1M2}
\sum_{m\leq M}\psi(- s^\gamma m^{2\gamma})
\ll \left| E_1(N,H,c)\right| + \left| E_2(N,H,c)\right| + M H^{-1},
\end{equation}
where
$$
E_1(N,H,c)= \sum_{m\leq M}\sum_{0<|h|\leq H}a(h)e(-hs^\gamma m^{2\gamma})
$$
and
$$
E_2(N,H,c)=  \sum_{m\leq M}\sum_{0<|h|\leq H}b(h)e(-hs^\gamma m^{2\gamma})
$$
(the term $M H^{-1}$ corresponds to the choice $h=0$ in  the summation on  the right hand side
of Lemma~\ref{lem:Vaal}).
We deal with $E_1(N,H,c)$ first. Switching the summation, we get
$$
E_1(N,H,c)\ll\sum_{0<|h|\leq H}\frac{1}{h}\left|\sum_{m\leq M}e(h s^\gamma m^{2\gamma})\right|.
$$
For the inner sum over $m$, we have
$$
\sum_{m\leq M}e(h s^\gamma m^{2\gamma})\ll \log N\max_{1\ll L \ll M}\left|\sum_{L<m\leq 2L}e(hs^\gamma m^{2\gamma})\right|.
$$
Using an exponent pair $(\kappa,\lambda)$, see~\cite{GrKol,Huxley,IwKow,Mont}, we obtain
$$
\sum_{L<m\leq 2L}e(h s^\gamma m^{2\gamma})\ll (h s^\gamma L^{2\gamma-1})^{\kappa} L^{\lambda}.
$$
Then
$$
\sum_{m\leq M}e(h s^\gamma m^{2\gamma})\LLN
h^{\kappa}s^{(\kappa-\lambda)/2} N^{c(2\gamma\kappa-\kappa+\lambda)/2} ,
$$
which yields
$H=s^{(\lambda-\kappa-1)/(2+2\kappa)}N^{c(1+\kappa-2\gamma\kappa-\lambda)/(2+2\kappa)}$
$$
E_1(N,H,c)\LLN H^{\kappa}s^{(\kappa-\lambda)/2}N^{c(2\gamma\kappa-\kappa+\lambda)/2}.
$$
By a similar argument, we can also get
$$
E_2(N,H,c)\LLN H^{\kappa}s^{(\kappa-\lambda)/2}N^{c(2\gamma\kappa-\kappa+\lambda)/2}.
$$
Applying Lemma~\ref{lem: Optim} to the bounds on terms in~\eqref{eq:M1M2},
 we obtain
$$
\sum_{m\leq s^{-1/2}N^{c/2}}\psi\(-s^\gamma m^{2\gamma}\)\LLN
s^{-\rho_1(\kappa,\lambda)}N^{\vartheta_1(c,\kappa,\lambda)}+s^{-\rho_2(\kappa,\lambda)}N^{\vartheta_2(c,\kappa,\lambda)}.
$$
Now the result follows from~\eqref{eq:Expand} and~\eqref{eq:MainTerm}.

\subsection{Proof of Theorem~\ref{thm: U-Bound} and~\ref{thm: U-Bound Aver}}

We fix some integer $P$ with $2\leq P \leq N$ (to be optimised later). It is also clear that
it is enough to obtain the desired bounds for $ Q_c^*(s;N)$, defined by~\eqref{eq:Q*}.

Using Corollary~\ref{cor:Sq Sieve PS} and then opening the square, changing the order
of summation and separating the diagonal terms (with
the total contribution at most $NP^{1+o(1)}$) , we obtain
\begin{equation*}
\begin{split}
Q_c^*(s; N) & \ll P^{-2}\sum _{\substack{P<\ell, p \leq 2P\\ \ell \neq p}}\log\ell \log p
\left|\sum_{N/2<n\leq N}
\(\frac{s \lfloor n^c\rfloor}{\ell p }\)\right|+  P^{-1}N^{1+o(1)}\\
&\LLN P^{-2}\sum _{\substack{P<\ell, p \leq 2P\\ \ell \neq p}} \left|\sum_{N/2<n\leq N}
\(\frac{ \lfloor n^c\rfloor}{\ell p }\)\right|+  P^{-1}N.
\end{split}
\end{equation*}
We remark that $s$ is no present anymore in the expression on the right hand side,
and thus the estimates below are uniform in $s$.

Note that the Jacobi symbols here are primitive characters thus~\eqref{eq:Betac}
applies and yields
 $$
Q_c^*(s; N) \LLN N^{1-\beta(c)}P+NP^{-1}.
$$
Taking $P=N^{\beta(c)/2}$, we get
$$
Q_c^*(s; N) \LLN N^{1-\beta(c)/2},
$$
which concludes the proof Theorem~\ref{thm: U-Bound}.

To prove Theorem~\ref{thm: U-Bound Aver} and we only need to consider
$$
\fQ_c^*(S,N)=\sum\limits_{\substack{s\leq S\\s~\text{square-free}}}Q_c^*(s; N).
$$
By Corollary~\ref{cor:Sq Sieve PS} again, we have
$$
\fQ_c^*(S,N) \LLN P^{-2}\sum _{\substack{P<\ell, p \leq 2P\\ \ell \neq p}}
\left|\sum_{\substack{s\leq S\\s~\text{square-free}}}\(\frac{s}{\ell p}\)\sum_{N/2<n\leq N}\(\frac{\lfloor n^c\rfloor}{\ell p }\)\right|+ P^{-1}SN.
$$
Applying~\eqref{eq:Betac}, we see that
$$
\begin{aligned}
\fQ_c^*(S,N)
&\LLN P^{-1}N^{1-\beta(c)}\sum _{\substack{r\leq 4P^2\\r~\text{square-free}}}
\left|\sum_{\substack{s\leq S\\s~\text{square-free}}}\(\frac{s}{r}\)\right|+P^{-1}SN
\end{aligned}
$$
Now by the Cauchy inequality, Lemma~\ref{Heath-Brown char sum} and choosing an optimal $P$, we get
$$
\fQ_c^*(S,N)\LLN (PS^{1/2}+S)N^{1-\beta(c)}+P^{-1}SN\ll SN^{1-\beta(c)}+S^{3/4}N^{1-\beta(c)/2},
$$
which yields Theorem~\ref{thm: U-Bound Aver}.
\subsection{Proof of Theorem~\ref{thm: U-Asymp Aver}}
\subsubsection{Preliminaries}
We proceed as in the proof of Theorem~\ref{thm: Asymp s}, then we have
\begin{equation}
\label{Start point}
\fQ_c(S,N) =S_0-E_1+E_2+O(1),
\end{equation}
where
$$
S_0=\sum\limits_{\substack{sm^2\leq N^c\\s\leq S\\s~\text{square-free}}}\left((sm^2+1)^{\gamma}-s^{\gamma}m^{2\gamma}\right),
$$
contributes to the main term and
$$
E_1=\sum\limits_{\substack{sm^2\leq N^c\\s\leq S\\s~\text{square-free}}}\psi(-s^\gamma m^{2\gamma})
\quad \text{and} \quad
E_2=\sum\limits_{\substack{sm^2\leq N^c\\s\leq S\\s~\text{square-free}}}\psi(-(sm^{2}+1)^\gamma).
$$
contribute to the error term.

\subsubsection{Evaluation of the main term $S_0$}
Using~\eqref{eq:MainTerm}, we compute $S_0$ directly as follows:
\begin{equation}
\label{S_0=}
S_0=\gamma(2\gamma-1)^{-1}N^{1-c/2} \Phi(S)+O(SN^{1-c}).
\end{equation}

\subsubsection{Reductions in the error terms  $E_1$ and $E_2$} By Lemma~\ref{lem:Vaal}, we obtain the following
analogue of~\eqref{eq:M1M2}:
\begin{equation}\label{S1=}
E_1\ll |E_{11}|+|E_{12}|+ H^{-1}N^{c/2}S^{1/2},
\end{equation}
where
$$
E_{11}=\sum\limits_{\substack{sm^2\leq N^c\\s\leq S\\s~\text{square-free}}}\sum\limits_{0<|h|\leq H}a(h)e(-hs^{\gamma}m^{2\gamma})
$$
and
$$
E_{12} = \sum\limits_{\substack{sm^2\leq N^c\\s\leq S\\s~\text{square-free}}}\sum\limits_{0<|h|\leq H}b(h)e(-hs^{\gamma}m^{2\gamma})
$$
for some $H\geq 2$. Using the same $H$, we  also have
\begin{equation}\label{S2=}
E_2 \ll |E_{21}|+|E_{22}|+ H^{-1}N^{c/2}S^{1/2},
\end{equation}
where
$$
E_{21}=\sum\limits_{\substack{sm^2\leq N^c\\s\leq S\\s~\text{square-free}}}\sum\limits_{0<|h|\leq H}a(h)e\left(-h(sm^2+1)^{\gamma}\right)
$$
and
$$
|E_{22}|= \sum\limits_{\substack{sm^2\leq N^c\\s\leq S\\s~\text{square-free}}}\sum\limits_{0<|h|\leq H}b(h)e\left(-h(sm^2+1)^{\gamma}\right)
$$
As usual, the sums $E_{12}$ and $E_{22}$ can be estimated similarly as $E_{11}$ and $E_{21}$, respectively, and by partial summation, $E_{21}$ can be converted to exponential sums which is similar to $E_{11}$
(see~\cite[Section~2]{HB1} for details).  In particular, we obtain same upper bounds for $E_{11}$, $E_{12}$, $E_{21}$
and $E_{22}$. Hence we only concentrate on the sum $E_{11}$.

\subsubsection{Estimating $E_{11}$}
Using
$$
\mu^2(s)=\sum\limits_{s=rd^2}\mu(d),
$$
we can write
$$
E_{11}=\sum\limits_{0<|h|\leq H}a(h)\sum_{\substack{rd^2m^2\leq N^c\\rd^2\leq S}}\mu(d)e(-hr^{\gamma}d^{2\gamma}m^{2\gamma}).
$$
Then, splitting the ranges of variables into dyadic ranges,  for  some real positive
parameters $R$, $D$ and $M$, satisfying
\begin{equation}\label{Conditions}
RD^2M^2\ll N^c \mand  RD^2\ll S,
\end{equation}
we obtain
\begin{equation}\label{E_11=}
E_{11}\LLN \sum_{0<|h|\leq H}h^{-1}\left|S(R,D,M;h)\right|,
\end{equation}
where
$$
S(R,D,M;h)=\sum\limits_{\substack{r\sim R,d\sim D,\ m\sim M,\\rd^2m^2\leq N^c,\ rd^2\leq S}}
\mu(d)e(-hr^{\gamma}d^{2\gamma}m^{2\gamma}).
$$

Now we estimate $S(R,D,M;h)$. Clearly we can assume that   $1<c<2$  as  for $c>2$ the result is
trivial (due to the presence of the term $S^{1/5}N^{(1+2c)/5} > N^{(1+2c)/5} \geq N$ for $c>2$).
We can remove the restrictive conditions
$$
rd^2m^2\leq N^c \mand rd^2\leq S
$$
at the cost of a small factor $(SN)^{o(1)}$ in a standard way (see,
for example,~\cite[Sections~2.3 and~3.2]{Harman}), which yields
$$
S(R,D,M;h)\LLN \left| \sum_{r\sim R,d\sim D,m\sim M}
\alpha_1(R)\alpha_2(D)\alpha_3(M)e(hr^{\gamma}d^{2\gamma}m^{2\gamma})\right|+M
$$
with some coefficients  $|\alpha_i(n)|\leq 1$ for $n\in\mathbb{N}$ and $i=1,2,3$.
Applying Lemma~\ref{Iwaniec's lemma} to the right hand side of the above formula, we get
\begin{align*}
S(R,D,M;h)\LLN&(hR^{\gamma}D^{2\gamma}M^{2\gamma})^{1/4}R^{1/2}(DM)^{3/4}+R^{7/10}DM\\
&\qquad +R(DM)^{3/4}+(hR^{\gamma}D^{2\gamma}M^{2\gamma})^{-1/4}R^{11/10}DM.
\end{align*}
Noting $\gamma>1/2$, it is easy to check that the fourth term can be absorbed by the third term on the right side. Thus by conditions~\eqref{Conditions}, we have
\begin{equation}\label{Bound2 for E_11}
E_{11}\LLN H^{1/4}S^{1/8}N^{1/4+3c/8}+S^{1/5}N^{c/2}+S^{5/8}N^{3c/8}.
\end{equation}

\subsubsection{Concluding the proof}
Bound~\eqref{Bound2 for E_11} together with~\eqref{S1=} and~\eqref{S2=}  yields
\begin{equation*}
\begin{split}
|E_{1}|+|E_2|\LLN &H^{1/4}S^{1/8}N^{1/4+3c/8}+H^{-1}N^{c/2}S^{1/2}\\
&+S^{1/5}N^{c/2}+S^{5/8}N^{3c/8}.
\end{split}
\end{equation*}
Now Lemma \ref{lem: Optim} gives
\begin{equation}
\label{eq: Bound E1E2}
|E_{1}|+|E_2|\LLN S^{1/5}N^{(1+2c)/5}+S^{5/8}N^{3c/8}+S^{1/8}N^{1/4+3c/8},
\end{equation}
where  the term $S^{1/5}N^{c/2}$ is absorbed by $S^{1/5}N^{(1+2c)/5}$, since we suppose $1<c<2$. Using the bound~\eqref{eq: Bound E1E2} together with ~\eqref{eq:Phi(S)}, ~\eqref{Start point} and~\eqref{S_0=}, and noting the contribution of $O(\log S)$ in ~\eqref{eq:Phi(S)} can also be absorbed by $S^{1/5}N^{(1+2c)/5}$, we obtain the desired result.
\section*{Acknowledgement}

The first two authors gratefully acknowledge the support, the hospitality
and the excellent conditions at the School of Mathematics and Statistics of UNSW during their visit.

This work was supported by NSFC Grant 11401329 (for K. Liu), by ARC Grant DP140100118 (for I. E. Shparlinski), the Natural Science Foundation of Shaanxi province of China Grant~2016JM1017(for T.~P.~Zhang).


\end{document}